\documentclass[11pt]{article}
\usepackage{amsfonts,amssymb,amsthm,amsmath}

\newtheorem*{theorem*}{Theorem}
\newtheorem*{MainTheor}{Main Theorem}
\newtheorem*{lemma}{Lemma}
\newtheorem*{lemma1}{Lemma 1}
\newtheorem*{lemma2}{Lemma 2}

\newtheorem*{conjecture1}{Conjecture 1}
\newtheorem*{conjecture2}{Conjecture 2}
\newtheorem*{intconj}{The Integral Conjecture}

\newtheorem*{corollary}{Corollary}
\theoremstyle{remark}
\newtheorem*{remark}{Remark}

\long\def\comment#1\endcomment{\relax}

\theoremstyle{definition}
\newtheorem*{definition}{Definition}

\DeclareMathOperator{\Dif}{Dif}
\DeclareMathOperator{\Hoch}{Hoch}
\DeclareMathOperator{\Cycl}{Cycl}

\DeclareMathOperator{\str}{str}
\DeclareMathOperator{\Tr}{Tr}
\DeclareMathOperator{\Id}{Id}

\newcommand{\mb}{{\bullet}}
\newcommand\End{\mathrm{End}}

\newcommand\E{\mathcal E}

\newcommand{\gl}{\mathfrak{gl}}
\newcommand{\cO}{\mathcal{O}}
\newcommand{\cA}{\mathcal{A}}
\newcommand{\cB}{\mathcal{B}}
\newcommand{\F}{\mathcal{F}}

\newcommand{\cD}{\mathcal{D}}

\newcommand{\Dolb}{\mathrm{Dolb}}

\newcommand{\K}{\mathrm{Ker}}
\newcommand{\CH}{\mathrm{H}}

\newcommand{\HH}{\mathrm{HH}}

\newcommand{\I}{\mathrm{Im}}

\def\wtilde#1{\widetilde{#1}\vphantom{#1}}

\title{ {\huge Riemann-Roch-Hirzebruch theorem and Topological Quantum Mechanics}}
\author{ Boris Feigin, Andrey Losev, and Boris Shoikhet} \date{}

\begin{document}\maketitle

\begin{abstract}
In the present paper we discuss an independent on the
Grothendieck-Sato isomorphism approach to the Riemann-Roch-Hirzebruch formula
for an arbitrary differential operator. Instead of the
Grothendieck-Sato isomorphism, we use the Topological Quantum Mechanics
(more or less equivalent to the well-known constructions with the Massey operations
from [KS], [P], [Me]). The statement that the Massey operations can "produce"
the integral in some set-up, has an independent from the RRH
theorem interest.

We finish the paper by some open questions arising when the main
construction is applied to the cyclic homology (instead of the Hochschild homology).
\end{abstract}

\section*{Introduction}
\subsection*{}
Let $X$ be an $n$-dimensional compact complex manifold, $E$ be a holomorphic
vector bundle over $X$. Consider the sheaf $\E$ of holomorphic sections of
the bundle $E$. Denote by $\Dif(\E)$ the algebra of global holomorphic
differential operators acting in the sheaf $\E$.

Each $\cD\in\Dif(\E)$ acts on $\CH^i(X,\E)$, $i\ge 1$. Define the super-trace
of $\cD$ as
\begin{equation}\label{eqr1.1}
\str(\cD):=\sum_{n=0}^\infty(-1)^i\Tr|_{\CH^i(X,\E)}(\cD)
\end{equation}
As each space $\CH^i(X,\E)$ is finite-dimensional because $X$ is compact, and
$\CH^i(X,\E)=0$ for $i>n$, the super-trace $\str(\cD)$ is well-defined.

It is clear that $\str\colon\Dif(\E)\to\mathbb{C}$ is a {\it trace } on the
algebra $\Dif(\E)$, that is, $\str([\cD_1,\cD_2])=0$ for any
$\cD_1,\cD_2\in\Dif(\E)$ where $[\cD_1,\cD_2]=\cD_1\cdot \cD_2-\cD_2\cdot \cD_1$.

When $\cD$ is the identity differential operator in the sheaf $\E$,
$\cD=\Id$, the super-trace (\ref{eqr1.1}) is equal to the Euler
characteristic $\chi(\E)=\sum_{i\ge 0}(-1)^i\dim \CH^i(X,\E)$
of the sheaf $\E$. Then the Riemann-Roch-Hirzebruch theorem expresses
$\chi(\E)$ via topological invariants of the manifold $X$ and of the bundle
$E$. Namely, $\chi(\E)=\int_X {Td}(T_X)\cdot{ch}(E)$ where ${Td}$ and ${ch}$
are the Todd genus and the Chern character, respectively.

In the present paper we find a formula for $\str(\cD)$ for an arbitrary
$\cD$ using
\begin{itemize}
\item[1)] a technical assumption that $\chi(\E)\ne 0$,
\item[2)] the Riemann-Roch-Hirzebruch theorem itself.
\end{itemize}
To be honest we
should warn the reader that the formula  for
$\str(\cD)$ involving explicitely the Hochschild cocycle from [FFSh] will not appear in this paper. Instead of it, we deal with
another formula, see (\ref{eqsem8}) below, which is equivalent to the formula
for $\str(\cD)$ from the Hochschild cocycle.

We decided to write this
paper by several reasons. At first, the proof here uses absolutely different
from the existing approaches
ideas and maybe even more suitable for further generalisations (for higher
algebraic K-theory, etc). At second, the integral here appears as a
consequence of a very general construction of the Topological Quantum
Mechanics ($\sim$ Massey operations in homotopical algebra). We have in mind
the existence of other applications of this scheme. The integral will be
replaced by something different, but from our point of view, will be given
by an analogous construction. This algebraic (physical) definition of the
integral and how it works in the framework of the RRH theorem was one of our
main motivations to write this paper.
\subsection*{The main idea}
Let us outline the main idea of our approach here.

As the super-trace $\str(\cD)$ is a trace on the associative algebra
$\Dif(\E)$, it is natural to consider it as a linear functional on 0-th
Hochschild homology $\HH_0(\Dif(\E))$. The Hocschild homology
$\HH_\mb(\Dif(\E))$ looks very big, and we have no way to compute it. In
particular, the associative algebra $\Dif(\E)$ can admit many other traces
besides $\str$. We want to do several things to improve the situation,
namely, to redefine the Hochschild homology $\HH_\mb(\Dif(\E))$ such that
they will become computable and nice. As usual, we have two ways to do that.
First, we can replace the algebra $\Dif(\E)$ to a natural dg algebra
$\Dif^\mb(\E)$ such that $\CH^0(\Dif^\mb(\E))\simeq \Dif(\E)$. After that,
we can hope that the Hochschild (and any other) homology
$\HH_\mb(\Dif^\mb(\E))$ is more computable. Second, we can replace the
concept of the Hochschild chain complex using some completions. In the
definition of the Hochschild homology we meant before, the Hochschild chain
of an algebra $A$ are elements from the (algebraic) tensor products
$A^{\otimes k}$. We can consider completed tensor products, and the homology
can change after that (and it is really the case).

In the first step, we set $\Dif^\mb(\E)$ be the Dolbeault complex of the
sheaf $\wtilde{\Dif(\E)}$ of holomorphic differential operators in the sheaf $\E$,
\begin{equation}\label{eqsem1}
\Dif^\mb(\E)=\Gamma_X(\Dolb^\mb(X, \wtilde{\Dif(\E)}))
\end{equation}
where
\begin{equation}\label{eqsem2}
\Dolb^\mb(X,
\wtilde{\Dif(\E)}):=\Dolb^\mb(X,\cO)\otimes_{\cO}\wtilde{\Dif(\E)}
\end{equation}
(Here $\Dolb^\mb(X,\cO)$ is the $\bar{\partial}$-resolution of the structure sheaf on the manifold
$X$).
Still now, the algebraic Hochschild homology $\HH_\mb(\Dif^\mb(\E))$
do not seem to be regular.

In Section 3.2 we prove the following theorem (see Theorem 3.2):
\begin{theorem*}
There exists a completion of all tensor products $\widehat{\Dif^\mb(\E)^{\otimes k}}$, $k\ge 1$, such that for the corresponding Hochschild homology
$\widehat{\HH_\mb}(\Dif^\mb(\E))$ one has:
\begin{equation}\label{eqsem3}
\widehat{\HH_{-k}}(\Dif^\mb(\E))=\CH^{2n-k}(X,\mathbb{C})
\end{equation}
where $n=\dim_{\mathbb{C}}(X)$. We use the convention (very natural from the
point of view of homological algebra) that the Hochschild chains of a usual
algebra $A$ (concentrated in degree 0) are negatively graded
$(\deg(A^{\otimes (k+1)})=-k$).
\end{theorem*}
We are interesting basically in traces, therefore, in 0-th Hochschild
homology.
Clearly we have a map
\begin{equation}\label{eqsem4}
\vartheta_1\colon \HH_0(\Dif^\mb(\E))\to \widehat{\HH_0}(\Dif^\mb(\E))
\end{equation}
because we have the inclusion of the algebraic tensor powers of an algebra
to completed tensor powers, and this map is compatible with the
differential.
Besides of it, we have the map
\begin{equation}\label{eqsem5}
\vartheta_0\colon\HH_0(\Dif(\E))\to \HH_0(\Dif^\mb(\E))
\end{equation}
induced by the imbedding of algebras $i\colon\Dif(\E)\hookrightarrow
\Dif^\mb(\E)$.
Denote by $\vartheta$ the composition
$\vartheta=\vartheta_1\circ\vartheta_0$.
Then we have the map
\begin{equation}\label{eqsem6}
\vartheta\colon\HH_0(\Dif(\E))\to\widehat{\HH_0}(\Dif^\mb(\E))\simeq\CH^{2n}(X)
\end{equation}
(the last isomorphism is stated by the Theorem above, the explicit form of
this isomorphism will be specified later).
As any global holomorphic differential operator $\cD\in\Dif(\E)$ defines an
element in $\HH_0(\Dif(\E))$, the map $\vartheta$ gives a map
\begin{equation}\label{eqsem7}
\cD\mapsto [\cD]\in\CH^{2n}(X)
\end{equation}
Our subject in this paper is to describe the last map, that is, to compute
$[\cD]$ for any holomorphic differential operator $\cD\in\Dif(\E)$.
\begin{conjecture1}
Let $X$ be a compact manifold. Then
\begin{equation}\label{eqsem8}
\int_{<X>}[\cD]=\str(\cD)
\end{equation}
\end{conjecture1}
Here we prove
the Conjecture using the technical assumption that
$\chi(\E)\ne 0$. In the proof here we use the Riemann-Roch-Hirzebruch
Theorem.
\begin{MainTheor}
The Conjecture above is true when $\chi(\E)\ne 0$.
\end{MainTheor}

Our main step in the proof of the Main Theorem is to lift the super-trace
$\str$ first to $\HH_0(\Dif^\mb(\E))$, and second to $\widehat{\HH_0}
(\Dif^\mb(\E))$. Denote the last linear functional by $\widehat{\str}$.
Then we have two linear functionals on $\CH^{2n}(X)$, the first is the integral,
and the second is $\widehat{\str}$. As $\CH^{2n}(X)$ is 1-dimensional, to
prove that they coincide it is enough to have an element $\cD\in\Dif(\E)$
such that the values of the both functionals on $\cD$ coincide and are {\it
non-zero}. We take $\cD=\Id$. Then we need to compute $[\Id]$. We do it in
[FFSh], this the only place where we refer to [FFSh]. We prove there that
$[\Id]={Td}(T_X)\cdot {ch}(E)$. Then we apply the RRH Theorem to prove that
the two linear functionals on $\CH^{2n}(X)$ coincide on $\cD=\Id$. When this
value is 0, we can not deduce that the two linear fuctionals coincide, it
just means that $\cD=\Id$ belongs to the kernel of the map $\vartheta$.
To be more precise, let us reformulate the Conjecture above:
\begin{conjecture2}
The map $\widehat{\str}\colon \CH^{2n}(X)\to\mathbb{C}$ coincides with the
integral map.
\end{conjecture2}
It is clear that the first Conjecture follows from the second one. Let us
notice that Conjecture 2 by its formulation is completely independent on
the RRH theorem. In our limping approach to this Conjecture we use the RRH
theorem to prove it in the case when $\chi(\E)\ne 0$. In the case
$\chi(\E)=0$ Conjecture 2 still remains to be open.

In this approach, the most difficulty (and actually the subject of this
paper) is to lift the super-trace $\str$ defined on $\HH_0(\Dif(\E))$ to a
map $\widehat{\str}$ defined on $\widehat{\HH_0}
(\Dif^\mb(\E))$. Actually this is the only difficulty. To do it, we should
"work on the level of complexes". We do it using the general formalism of
Topological Quantum Mechanics, introduced by one of us (A.S.L) with
Vyacheslav Lysov [Lys]. This formalism is a nice way to reformulate the well-known
construction with the Massey operations (see [KS], [P], [Me]).
The statement that the Massey operations can give the analytic integral,
looks to be new (the statement of Conjecture 2 we prove in the particular
case when $\chi(\E)\ne 0$). It is very interesting to find possible
generalizations of this construction, as well as to find a more direct
(independent on the RRH Theorem) approach to the Conjecture 2.
\newpage

\section{Hochschild homology and super-traces}
Recall that the Hochschild homology of an associative algebra $A$ is the
(co)homology of the Hochschild chain complex $\Hoch_\mb(A)$. This is the
complex
$$
\dots\rightarrow
\Hoch_2(A)\overset{d_{\Hoch}}\rightarrow\Hoch_1(A)\overset{d_{\Hoch}}\rightarrow
\Hoch_0(A)\rightarrow 0
$$
where $\Hoch_k(A)=A^{\otimes (k+1)}$, and the differential $d_{\Hoch}\colon
A^{\otimes (k+1)}\to A^{\otimes k}$ is defined as follows:
\begin{equation}\label{eqr1.4}
\begin{aligned}
d&_{\Hoch}(a_0\otimes a_1\otimes\dots\otimes a_k)=\\
&a_0a_1\otimes a_2\otimes\dots\otimes a_k
-a_0\otimes (a_1a_2)\otimes\dots\otimes a_k
+\dots\dots\dots\\
&+(-1)^{k-1}a_0\otimes\dots\otimes (a_{k-1}a_k)
+(-1)^k a_ka_0\otimes\dots\otimes a_{k-1}
\end{aligned}
\end{equation}
Denote by $\HH_{-i}(A)$ the $i$-th Hochschild homology of the algebra $A$.
In particular, $\HH_0(A)=A/[A,A]$. Therefore, a trace on the algebra $A$ is
the same that an element from the dual space $(\HH_0(A))^*$. Therefore, the super-trace
$\str\colon \Dif(\E)\to\mathbb{C}$ constructed
above is a canonical element in $(\HH_0(\Dif(\E))^*$.

We want to consider the Hochschild homological differential $d_{\Hoch}$ as a
differential of degree $+1$. For this, we use the convention that $\deg
A^{\otimes (k+1)}=-k$.

We have the definition of the Hochschild homological complex also in the
case of a differential graded associative algebra $A^\mb$. It is a bicomplex
with the two differentials $d_{A^\mb}$ and $d_{\Hoch}$. So, if $A^\mb$ is
$\mathbb{Z}_{\ge 0}$-graded, the bicomplex $\Hoch_\mb(A^\mb)$ is non-zero in
the quarter-plane $\{x\le 0,\ y\ge 0\}$.

\begin{definition}
Let $A^\mb$ be a dg associative algebra. Then a trace on $A^\mb$ is a linear
functional $\Tr\colon\mathbb{H}^0(\Hoch_\mb(A^\mb))\to\mathbb{C}$.
\end{definition}
(We denote by $\mathbb{H}^\mb$ the hyper-cohomology of a bicomplex, which
is, by definition, the cohomology of the corresponding ordinary complex).

Our first goal in this paper is to construct a trace on the dg algebra $\Dif^\mb(\E)$
such that its restriction to $\Dif(\E)\subset \Dif^\mb(\E)$ is the
super-trace $\str$
constructed above.

It means that we want to construct linear functionals
$$
\Theta_k\colon[(\Dif^\mb(\E))^{\otimes k}]^{k-1}\to\mathbb{C}, k\ge 1
$$
(here $[a\otimes\dots\otimes z]^i$ means that $\deg a+\dots+\deg z=i$)
such that for any $\alpha\in [(\Dif^\mb(\E))^{\otimes k}]^{k-2}$ one has
\begin{equation}\label{eqr1.7}
\Theta_{k-1}(d_{\Hoch}\alpha)+\Theta_k(\bar{\partial}\alpha)=0,\ \ \ k\ge 2
\end{equation}

In particular, in the simplest case $k=2$ one has:
\begin{equation}\label{eqr1.7}
\Theta_1([\cD_1,\cD_2])+\Theta_2(\bar{\partial}\cD_1\otimes
\cD_2)+\Theta_2(\cD_1\otimes \bar{\partial}\cD_2)=0
\end{equation}
where $\cD_1, \cD_2$ are $C^{\infty}$-differential operators actiong in the
sheaf $\E$. When $\cD_1, \cD_2$ are holomorphic, $\bar{\partial}\cD_i=0$, we
obtain the usual trace condition $\Theta_1([\cD_1,\cD_2])=0$.

Denote by $K^\mb=\Gamma_X(\Dolb^\mb(X,\E))$ the Dolbeault complex of the sheaf
$\E$.
Let us suppose we have a decomposition of {\it complexes }
$K^\mb=K^\mb_0\bigoplus K^\mb_1$ where $K^\mb_0\simeq H^\mb(K^\mb)$,
$d_{K^\mb_0}=0$, and $H^\mb(K^\mb_1)=0$.
For any $C^{\infty}$ differential operator $\cD$ define $\Theta_1(\cD)$ as
\begin{equation}\label{eqr1.8}
\Theta_1(\cD):=\str(\Pi_{K^\mb_0}\cD\Pi_{K^\mb_0})
\end{equation}
where $\Pi_{K^\mb_0}, \ \Pi_{K^\mb_1}$ are the projectors of $K^\mb$ to
$K^\mb_0,\ K^\mb_1$, correspondingly, and the super-trace here is the
super-trace of the endomorphism of $K^\mb_0$.

It is clear that for a holomorphic $\cD$, $\Theta_1(\cD)$ coincides with the
usual super-trace. For non-holomorphic differential operators,
$\Theta_1([\cD_1,\cD_2])\neq 0$, and the problem is to find $\Theta_2$ such that
equation (\ref{eqr1.7}) holds.

There is a  theory due to Andrey Losev and Vyacheslav Lysov [Lys],
called topological quantum mechanics, which allows us to construct
all $\Theta_k,\ k\geq 1$.

\section{The Topological Quantum Mechanics}\label{sec2}
\subsection{The general set-up}\label{subsec2.1}
Suppose $\cA^\mb$ is a dg associative algebra, and $M^\mb$ is a dg module
over it. It means in particular that there is a map
$$
\rho\colon\cA^\mb\otimes M^\mb\to M^\mb
$$
which is a map of {\it complexes}. Denote by $d_{\cA}$ the differential in
$\cA^\mb$, and by $Q$ the differential in $M^\mb$. Denote by $\rho(a)$, $a\in \cA^\mb$,
the corresponding endomorphism of $M^\mb$.
Then
\begin{equation}\label{eqr1.9}
\rho(d_{\cA}(a))=\{ Q,\rho(a)\}
\end{equation}
Suppose also that $M^\mb$ has {\it finite-dimensional cohomology}.

Our main example is the case $\cA^\mb=\Dif^\mb(\E),\ M^\mb=K^\mb$ in the notations
of previous Sections.

Choose a decomposition $M^\mb=M^\mb_0\oplus M^\mb_1$ where $M^\mb_0\simeq
H^\mb(M^\mb)$, $d_{M^\mb_0}=0$, and $H^\mb(M^\mb_1)=0$.
The topological quantum mechanics constructs an $A_\infty$-morphism
$\F\colon\cA^\mb\to\End M^\mb_0$. That is, it constructs a collection of maps
\begin{equation}\label{eqsem10}
\F_k\colon\cA^{\mb\otimes k}\to \End M_0^\mb[1-k],\ \ k\ge 1
\end{equation}
satisfying the following relations:
\begin{equation}\label{eqsem11}
\begin{aligned}
\ &[\F_{k-1}((a_1\cdot a_2)\otimes a_3\otimes\dots\otimes a_k)\mp
\F_{k-1}(a_1\otimes (a_2\cdot a_3)\otimes\dots\otimes
a_k)\pm\dots\\
& \pm\F_{k-1}(a_1\otimes  a_2\otimes\dots\otimes (a_{k-1}\cdot a_k))]\pm\\
& \pm\F_k(\bar{\partial}(a_1\otimes a_2\otimes\dots\otimes a_k))\pm\\
& \pm\sum_{\ell =1}^{k-1} \F_{\ell}(a_1\otimes\dots\otimes a_{\ell})\circ
\F_{k-{\ell}}(a_{\ell +1}\otimes\dots\otimes a_k)=0,\\
& k\ge 1
\end{aligned}
\end{equation}
Here in the last term $\circ$ denotes the composition in
$\End(M^\mb_0)$.
We construct such an $A_\infty$-morphism below in this Section.
\begin{remark}
In [KS], [P], [Me] the authors consider the following situation. Let
$\mathbf{B}^\mb$ be an associative dg algebra, and suppose we have a
decomposition $\mathbf{B}^\mb=\mathbf{B}^\mb_0\oplus\mathbf{B}^\mb_1$ where
$\CH^\mb(\mathbf{B}_0^\mb)\simeq \CH^\mb(\mathbf{B}^\mb)$ and
$\mathbf{B}^\mb_1$ is acyclic. Notice that we do {\it not} suppose that the
differential in $\mathbf{B}_0^\mb$ is zero and do {\it not} suppose that
$\mathbf{B}^\mb_0$ is a subalgebra. Then in the papers cited above the
authors construct an $A_\infty$-algebra structure on $\mathbf{B}^\mb_0$ in
these assumptions such that this $A_\infty$-algebra is $A_\infty$-isomorphic
to $\mathbf{B}^\mb$. It is more or less a construction of higher Massey operations. The relation with our situation is as follows: set
$\mathbf{B}^\mb=\cA^\mb\oplus M^\mb[-1]$ where the product on $\cA^\mb$ is
the initial product, the product of any two elements in $M^\mb[-1]$ is zero,
and the product of an element in $\cA^\mb$ with an element in $M^\mb[-1]$
is the module action. Then $\mathbf{B}^\mb$ is an associative dg algebra.
Set $\mathbf{B}_0^\mb=\cA^\mb\oplus M_0^\mb[-1]$ and $\mathbf{B}^\mb_1=M_1^\mb[-1]$.  It is clear that we are in the assumptions above. Then in the
papers cited above the authors construct an $A_\infty$-algebra structure on
$\cA^\mb\oplus M^\mb_0[-1]$ such that the product restricted to $\cA^\mb$ is the
initial product, and the higher products containing more than 1 element from
$M_0^\mb$ are 0. It is clear that this structure is exactly the same that an
$A_\infty$-morphism $\F\colon\cA\to\End M_0^\mb$. We are going to present an
exposition of this construction independent on the papers cited above.
\end{remark}
\bigskip
\bigskip
\bigskip
\bigskip
\bigskip

\subsection{Construction of $\F_k$}\label{sec2.2}
Let $a_1,\dots,a_k\in\cA^\mb$, $\tau_1 <\tau_2 <\dots <\tau_{k}$ be real
numbers. Choose a homotopy $\kappa\colon M^\mb\to M_1^\mb[-1]$  such that
$\{Q, \kappa\}=\Pi_{M_1^\mb}$ where $\Pi_{M_i^\mb}$ is the projection to $M_i^\mb$. Define the spaces
$$
{Conf}_k=\{\tau_1<\dots<\tau_k,\ \ \tau_i\in\mathbb{R}\}
$$
and
$$
C_k={Conf}_k/G^{(1)}
$$
where $G^{(1)}$ is the 1-dimensional group of shifts
$\tau_i\mapsto\tau_i+c,\ i=1,\dots,k$.

Set
\begin{multline}\label{eqr1.12}
\Omega_{a_1,\dots,a_k}=\Pi_{M^\mb_0}\circ \rho(a_k)\circ \exp[-d(\tau_{k}-\tau_{k-1})\kappa-
(\tau_{k}-\tau_{k-1})\Pi_{M^\mb_1}]\circ \rho(a_{k-1})\circ \\
\exp[-d(\tau_{k-1}-\tau_{k-2})\kappa-
(\tau_{k-1}-\tau_{k-2})\Pi_{M^\mb_1}]\circ \rho(a_{k-2})\circ \\
\dots\circ\exp[-d(\tau_2-\tau_1)\kappa-
(\tau_2-\tau_1)\Pi_{M^\mb_1}]\circ \rho(a_1)\circ \Pi_{M^\mb_0}\in \End^\mb(M^\mb_0)\otimes
\Omega^\mb(C_{k+1})\\
(a_1,\dots , a_k\in\cA^\mb)
\end{multline}
which is a non-homogeneous differential form with the top component of degree $k-1$ on the space
$C_{k}$ with values in $\End(M^\mb_0)$.

Denote $d_{\cA}\Omega_{a_1,\dots,a_k}=\sum_{i=1}^k\Omega_{a_1,
\dots ,d_{\cA}a_i,\dots,a_k}$. Denote by $d_\tau$ the de Rham differential
in the space $C_{k}$.
\begin{lemma}
\begin{equation}\label{eqrl1}
(d_{\cA}-d_{\tau})\Omega_{a_1,\dots,a_k}=0
\end{equation}
\begin{proof}
We use the identity for $t\in\mathbb{R}$
\begin{equation}\label{eqr1.13}
dt\kappa+t\Pi_{M^\mb_1}=\{Q+d_t,t\kappa\}
\end{equation}
In particular,
\begin{equation}\label{eqr1.14}
\{Q+d_t,\exp [dt\kappa+t\Pi_{M^\mb_1}]\}=0
\end{equation}
We have:
\begin{multline}\label{eqr1.15}
d_{\cA}\Omega_{a_1,\dots,a_k}=\\
\Pi_{M^\mb_0}\circ \rho(d_{\cA}(a_k))\circ \exp[-d(\tau_{k}-\tau_{k-1})\kappa-
(\tau_{k}-\tau_{k-1})\Pi_{M^\mb_1}]\circ \rho (a_{k-1})\circ \\
\exp[-d(\tau_{k-1}-\tau_{k-2})\kappa-
(\tau_{k-1}-\tau_{k-2})\Pi_{M^\mb_1}]\circ \rho(a_{k-2})\circ \\
\dots\circ\exp[-d(\tau_2-\tau_1)\kappa-
(\tau_2-\tau_1)\Pi_{M^\mb_1}]\circ \rho(a_1)\circ \Pi_{M^\mb_0}+\dots=\\
\bigskip
\Pi_{M^\mb_0}\circ \{Q,\rho(a_k)\}\circ \exp[-d(\tau_{k}-\tau_{k-1})\kappa-
(\tau_{k}-\tau_{k-1})\Pi_{M^\mb_1}]\circ \rho(a_{k-1})\circ \\
\exp[-d(\tau_{k-1}-\tau_{k-2})\kappa-
(\tau_{k-1}-\tau_{k-2})\Pi_{M^\mb_1}]\circ \rho(a_{k-2})\circ \\
\dots\circ\exp[-d(\tau_2-\tau_1)\kappa-
(\tau_2-\tau_1)\Pi_{M^\mb_1}]\circ \rho(a_1)\circ \Pi_{M^\mb_0}+\dots
\end{multline}
by (\ref{eqr1.9}).

Next, the r.h.s. of (\ref{eqr1.15}) is
\begin{multline}\label{eqr1.16}
-\Pi_{M_0}\circ \rho(a_k)\circ\{Q, \exp[-d(\tau_{k}-\tau_{k-1})\kappa-
(\tau_{k}-\tau_{k-1})\Pi_{M^\mb_1}]\}\circ \rho(a_{k-1})\circ \\
\exp[-d(\tau_{k-1}-\tau_{k-2})\kappa-
(\tau_{k-1}-\tau_{k-2})\Pi_{M^\mb_1}]\circ \rho(a_{k-2})\circ \\
\dots\circ\exp[-d(\tau_2-\tau_1)\kappa-
(\tau_2-\tau_1)\Pi_{M^\mb_1}]\circ \rho(a_1)\circ \Pi_{M^\mb_0}+\dots
\end{multline}
because $Q$ acts by 0 on $\End(M^\mb_0)$ and commutes with $\Pi_{M^\mb_0}$.
Now the claim of Lemma follows from (\ref{eqr1.14}).
\end{proof}
\end{lemma}

This Lemma means that the differential form $\Omega_{a_1,\dots,a_k}$ is "almost closed",
and we can apply to it the Stokes formula as usually in topological field
theories to get some algebraic identities (see, for instance, [K]).

We set
\begin{equation}\label{eqr1.17}
\F_k(a_1\otimes\dots\otimes a_k):=\int_{C_k}\Omega_{a_1,\dots,a_k}\in \End(M^\mb_0)
\end{equation}

The space $C_k$ is a smooth manifold of dimension $k-1$, and only the top
degree component (of degree $k-1$) of the non-homogeneous differential form
$\Omega_{a_1,\dots,a_k}$ contributes. A priori this integral could diverge
because the space $C_k$ is {\it not } compact. However, it will be clear
from the sequel, that the integral absolutely converges.

\subsection{The Stokes formula}\label{sec2.3}
Here we prove the following result:
\begin{theorem*}
The maps $\F_k\colon{\cA}^{\otimes k}\to\End(M^\mb_0)[1-k]$, $k\ge
1$, satisfy the followning relation:
\begin{equation}\label{eqr1.18}
\begin{aligned}
\ &[\F_{k-1}((a_1\cdot a_2)\otimes a_3\otimes\dots\otimes a_k)\mp
\F_{k-1}(a_1\otimes (a_2\cdot a_3)\otimes\dots\otimes
a_k)\pm\dots\\
& \pm\F_{k-1}(a_1\otimes  a_2\otimes\dots\otimes (a_{k-1}\cdot a_k))]\pm\\
& \pm\F_k(\bar{\partial}(a_1\otimes a_2\otimes\dots\otimes a_k))\pm\\
& \pm\sum_{\ell =1}^{k-1} \F_{\ell}(a_1\otimes\dots\otimes a_{\ell})\circ
\F_{k-{\ell}}(a_{\ell +1}\otimes\dots\otimes a_k)=0,\\
& k\ge 1
\end{aligned}
\end{equation}
$\mathrm{(}$Here in the last term $\circ$ denotes the composition in
$\End(M^\mb_0)\mathrm{)}$.

\begin{proof}
First, we consider a compactification $\overline{C_k}$ of the manifold
$C_k$ for any $k\ge 1$. This compactification is a manifold with corners
stratified by locally-closed strata isomorphic to products $C_{i_1}\times
C_{i_2}\times\dots\times C_{i_{\ell}}$. Let us first describe this
compactification informally. When two points in ${Conf}_k$ move close to
each other, we want to have a stratum of codimension 1. When $\ell$ points
move close to each other we want to have a stratum of codimension $\ell -1$.
Next, when the distance between two neihbour points tends to $\infty$, say
the distance between the $\ell$'th and the $(\ell+1)$'st points, then we want
to have the stratum $C_{\ell}\times C_{k-\ell}$. Next, for higher
degenerations the strata are labeled by trees, and so on.

There exists a nice explicit construction of $\overline{C_k}$, the {\it cube }.
Consider the cube generated by $k-1$ ortogonal vectors of length 1. The $i$'th vector has the length
$\tau_{i+1}-\tau_i$. We compactify the ray $\{\tau_{i+1}-\tau_i| \ \tau_{i+1}>\tau_i\}$ by a point at the infinity.
Then
it is clear that this compactification obeys all properties listed above.

It is clear that the differential form $\Omega_{a_1,\dots,a_k}$ constructed
above can be extended to a smooth closed differential form on the cube
$\overline{C_k}$. This proves, in particular, that the integrals
$\int_{C_k}\Omega_{a_1,\dots,a_k}$ absolutely converge.

Now we apply the Stokes formula to get (\ref{eqr1.18}).
We have from (\ref{eqrl1})
\begin{equation}\label{eqr1.19}
0=
\int_{\partial\overline{C_{k}}}\Omega_{a_1,\dots,a_k}-
\int_{\overline{C_k}}d_{\cA}\Omega_{a_1,\dots,a_k}
\end{equation}

Recall that $\Omega_{a_1,\dots,a_k}$  is a non-homogeneous form with the top
component of degree $k-1$. We consider in details now  the first summand in
the r.h.s. of (\ref{eqr1.19}). We can restrict ourselves only in homogeneous
component of degrees $k-1$ and $k-2$ in $\Omega_{a_1,\dots,a_k}$. In
the first summand in the r.h.s of (\ref{eqr1.19}) only the boundary strata
of codimension 1 in $\partial\overline{C_k}$ contributes. These boundary
strata are exactly the faces of codimension 1 of the {\it cube }.

The boundary strata corresponding to these faces are of the two types, S1 and S2
below.

S1) Two points in ${Conf}_k$ move close to each other. The boundary stratum
is $C_{k-1}$.

S2) The points are divided by two groups, $\{p_1,\dots,p_{\ell}\}$ and
$\{p_{\ell+1},\dots,p_k\}$, and the distance between the groups thends to infinity. The boundary stratum is $C_{\ell}\times
C_{k-\ell}$.

Geometrically, let the {\it cube} we consider be the unit cube in
$\mathbb{R}^{k-1}$. Then it has two "distignuished" vertices, namely,
$A_0=(0,0,\dots,0)$ and $A_\infty=(1,1,\dots, 1)$. (Here the point 1 on each
ray $\{\tau_{i+1}-\tau_i|\tau_{i+1}>\tau_i\}$ is the image of the "point"
$\tau_{i+1}-\tau_i=\infty$ under the compactification we consider). Then
each face of codimension 1 of the {\it cube} contains either the point $A_0$
or the
point $A_\infty$. The case S1) above describes exactly the faces which
contain the point $A_0$, while the case S2) describes exactly the faces
which contain the point $A_\infty$.

It is clear that the other degenerations have higher codimension and are
corresponded to the faces of the cube of higher codimension.

Now we want to compute the integral over $\partial\overline{C_k}$. We have:

\begin{multline}\label{eqr1.20}
\int_{\partial\overline{C_{k}}}\Omega_{a_1,\dots,a_k}=
\int_{\partial_{S1}\overline{C_{k}}}\Omega_{a_1,\dots,a_k}+
\int_{\partial_{S2}\overline{C_{k}}}\Omega_{a_1,\dots,a_k}
\end{multline}
We need to compute the restrictions $\Omega_{a_1,\dots,
a_k}|_{\partial_{S1}}$ and $\Omega_{a_1,\dots,
a_k}|_{\partial_{S2}}$.
\subsubsection{}\label{sssectr1.3.2.1}
\begin{lemma}
\begin{itemize}
\item[]
\item[$(i)$]$\Omega_{a_1,\dots,a_k}|_{\partial_{S1,i}}=\pm \Omega_{a_1,\dots,
a_ia_{i+1},\dots,a_k}$ as a differential form on the space $C_{k-1}$. Here we denote by $S1,i$ the boundary stratum of the
type S1 corresponding to the case when the $i$-th point moves close to the
$i+1$-st, $1\le i\le k-1$.
\item[$(ii)$]$\Omega_{a_1,\dots,a_k}|_{\partial_{S2,\ell}}=\pm\Omega_{a_1,\dots,a_\ell}
\boxtimes\Omega_{a_{\ell+1},\dots,a_k}$ as a differential form on the space
$C_{\ell}\times C_{k-\ell}$. Here $\boxtimes$ denotes the external product
of forms on the product of spaces, and the $\bf{composition}$ of the forms
with values in $\End M^\mb_0$.
\end{itemize}
\begin{proof}
We prove (i) and (ii) separately.

For the proof of (i) denote $t_i=\tau_{i+1}-\tau_i$, $i\le k-1$. We should
analyze the quantity $\exp(-dt_i\kappa-t_i\Pi_{M_1^\mb})$ when
$t_i\rightarrow 0$. Then $dt_i$ is identically 0, and the whole exponent is
equal to 1. Now the claim (i) of Lemma follows from the identity
$\rho(a_{i+1})\cdot\rho(a_i)=\rho(a_{i+1}\cdot a_i)$ (which holds because $M^\mb$ is a dg
module over the dg algebra $\cA^\mb$).

For the proof of (ii), we should consider the case when
$t_i\rightarrow\infty$. Then $dt_i$ is also 0 identically, because $t_i=\infty$
is a point of the comtified space. Now the claim follows from the identity
\begin{equation}\label{eqsemsuper}
\exp(-\infty\Pi_{M_1^\mb})=\Pi_{M_0^\mb}
\end{equation}
It is clear: the exponent maps to 0 any element of $M_1^\mb$, and for $m\in M^\mb_0$
all summands in the power series expansion of the exponent
$\exp(-\infty\Pi_{M_1^\mb})(m)$ are 0 except for the first which is the
identity operator. The equation (\ref{eqsemsuper}), and the Lemma are
proven.
\end{proof}
\end{lemma}
Theorem 2.3 is proven.
\end{proof}
\end{theorem*}
We have proved that the collection of maps $\{\F_k,\ k\ge 1\}$ define an
$A_\infty$-morphism $\F\colon\cA^\mb\to\End M^\mb_0$.
\begin{remark}
Notice that we never used in the course of the proof that $\Pi_{M_1}$ is the
projector. Namely, we can axiomatize the situation as follows. We have a
decomposition $M^\mb=M_0^\mb\oplus M_1^\mb$, $Q$ is the differential in
$M^\mb$, and both $M_0^\mb$, $M^\mb_1$ are subcomplexes. Let $\Pi_{M_0}$ be
the projector to $M_0$. Choose a homotopy $\kappa\colon M^\mb\to
M^\mb_1[-1]$ such that $\{Q,\kappa\}=\Delta$ where $\Delta\colon M^\mb\to
M^\mb$ is an operator such that:
\begin{itemize}
\item[(i)] $\Delta|_{M_0}=0$,
\item[(ii)] $\Delta(M_1)\subset M_1$,
\item[(iii)] $\exp (-t\Delta)|_{M_1}=0$ when $t\rightarrow\infty$.
\end{itemize}
The conditions (i)-(iii) together mean that the exponent $\exp
(-\infty\Delta)$ is the projector $\Pi_{M_0}$. In other words, we can
replace the projector $\Pi_{M_1}$ by the operator $\Delta$ above with the
weaker properties. Namely, we replace the property $\Pi_{M_1}|_{M_1}=\Id$ of
the projector by the weaker property (iii) of the operator $\Delta$. We use this Remark later in
Section 3.1.2.
\end{remark}

\subsection{Applications to Hochschild and cyclic homology, and to
supertraces}\label{sec2.4}
We prove here the following statement:
\begin{lemma}
Any $A_\infty$-morphism of associative dg algebras algebras $\F\colon\cA^\mb\to\cB^\mb$ induces a
map of the Hochschild complexes $\F_{\Hoch}\colon\Hoch_\mb(\cA^\mb)\to
\Hoch_\mb(\cB^\mb)$, and a map of the cyclic complexes $\F_{\Cycl}\colon
\Cycl_\mb(\cA^\mb)\to\Cycl_\mb(\cB^\mb)$.
\begin{proof}
First construct the map $\F_{\Hoch}\colon\Hoch_\mb(\cA^\mb)\to
\Hoch_\mb(\cB^\mb)$. Let $\Omega\in\cA^{\mb\otimes k}$. Define the set $\Xi_k$
of all ordered partitions $\xi$ of $k$, $\xi=(k_1,k_2,\dots, k_{\ell}),\
k=k_1+k_2+\dots+k_{\ell}$.
Suppose $\Omega=a_1\otimes a_2\otimes\dots\otimes a_k$. Define
\begin{equation}\label{eqsem20}
\begin{aligned}
\ &\F_{\Hoch}(a_1\otimes\dots\otimes
a_k)=\sum_{\xi\in\Xi_k}\{ [\F_{k_1}(a_1\otimes\dots\otimes a_{k_1})\otimes\\
&\otimes\F_{k_2}(a_{k_1+1}\otimes\dots\otimes a_{k_2})\otimes\dots\otimes
\F_{k_\ell}(a_{k_1+\dots +k_{\ell-1}+1}\otimes\dots\otimes a_k)]\\
&+\sum_{s=2}^{k_\ell}(-1)^{n(\xi, s,\Omega)}[\F_{k_\ell}(a_{k-k_\ell+s}\otimes
\dots\otimes a_k\otimes a_1\otimes\dots\otimes
a_{s-1})\otimes\\
&\otimes\F_{k_1}(a_s\otimes\dots\otimes a_{s+k_1-1})\otimes\\
&\otimes\dots\otimes \\
&\otimes \F_{k_{\ell-1}}(a_{s+k_1+\dots+k_{\ell-2}}\otimes\dots\otimes
a_{s+k_1+\dots+k_{\ell-2}+k_{\ell-1}-1})]\}
\end{aligned}
\end{equation}
In the formula above $\F_{k_1},\dots,\F_{k_\ell}$ are the components of the
$A_\infty$-morphism $\F$, and the sign $(-1)^{n(\xi, s, \Omega)}$ is defined
as follows. If $\deg a_i=0$ for all $i=1,\dots, k$, $n(\xi, s,
\Omega)=(k-k_\ell+s-1)(k_\ell-s+1)$, that is, the sign is equal to the sign of the
permutation $(1,2,\dots, k)\mapsto (k-k_\ell+s,\dots, k, 1,
2,\dots,k-k_\ell+s-1)$.
In the general case, we should take in the account the degrees of $a_i$'s.

The reader can easily check that the map $\F_{\Hoch}$ defines in fact a map
of the Hochschild complexes. Moreover, it also defines a map of the cyclic
complexes, so we can set $\F_{\Cycl}=\F_{\Hoch}$.
\end{proof}
\end{lemma}

Next, if we have any trace
$\Upsilon\colon\mathbb{H}^0(\Hoch_\mb(\cB^\mb))\to\mathbb{C}$, then the composition
$\Upsilon_{new}=\Upsilon\circ\F_{\Hoch}$ gives us a trace
$\mathbb{H}^0(\Hoch_\mb(\cA^\mb))\to\mathbb{C}$.

\section{Applications to the algebra $\cA^\mb=\Dif^\mb(\E)$}\label{sec3}
\subsection{The trace on the algebra $\Dif^\mb(\E)$ before the completion}
We do not know what is the algebraic (not completed) Hochschild homology
$\mathbb{H}^\mb(\Hoch_\mb(\Dif^\mb(\E)))$. Nevertheless, the theory developed in
Section 2 allows to construct a trace
$\str^\mb\colon\mathbb{H}^0(\Hoch_\mb(\Dif^\mb(\E)))\to\mathbb{C}$. Let us
describe this trace (in more details than before).

At first, we have an $A_\infty$-morphism $\F\colon\Dif^\mb(\E)\to
\End(K_0^\mb)$ where $K^\mb=\Gamma_X(\Dolb^\mb(X,\E))$. In our assumptions
that $X$ is compact the cohomology $\CH^\mb(K^\mb)=K^\mb_0$ is finite-dimensional.
We have the corresponding map
\begin{equation}\label{eqsem30}
\F_{\Hoch}\colon\Hoch_\mb(\Dif^\mb(\E))\to\Hoch_\mb(\End K_0^\mb)
\end{equation}
(see Section 2.4).
In particular, we have the map
\begin{equation}\label{eqsem31}
\F^0_{\Hoch}\colon\mathbb{H}^0(\Hoch_\mb(\Dif^\mb(\E)))\to\mathbb{H}^0(\Hoch_\mb(\End
K_0^\mb))
\end{equation}
Therefore, to construct the super-trace
$\str^\mb\colon\mathbb{H}^0(\Hoch_\mb(\Dif^\mb(\E)))\to\mathbb{C}$,
we need to construct a trace $\mathbb{H}^0(\Hoch_\mb(\End
K_0^\mb))\to\mathbb{C}$.

As $K^\mb_0=\CH^\mb K^\mb$ are finite-dimensional, we have a very simple situation:
the Hochschild homology of the algebra $\mathbf{A}^\mb$ of endomorphisms of
a $\mathbb{Z}$-graded {\it finite-dimensional} vector spece $V^\mb$. We can use
the following statement:
\begin{lemma}
In the assumptions above, the Hochschild homology
$\mathbb{H}^i(\Hoch_\mb(\mathbf{A}^\mb))$ is $\mathbb{C}$ for $i=0$, and is 0 for
$i\ne 0$. We can describe a canonical (nonzero) super-trace
$\str_{can}\colon\mathbb{H}^0(\Hoch_\mb(\mathbf{A}^\mb))\to\mathbb{C}$ as
follows: it is nonzero on $[\mathbf{A}^{\mb\otimes k}]^{k-1}$ only for
$k=1$, and for an endomorphism of degree 0 $m\in\mathbf{A}^0$, the
super-trace is
\begin{equation}\label{eqsem32}
\str_{can}(m)=\Tr|_{V^{even}}(m)-\Tr|_{V^{odd}}(m)
\end{equation}
\end{lemma}
It is a standard fact.
\qed

Now we want to pull-back the super-trace $\str_{can}$ by the $A_\infty$-morphism
$\F$. This pull-back was described in Section 2.4 above. The only
specification here is that the trace $\str_{can}$ we start with is nonzero
on $[\mathbf{A}^{\mb\otimes \ell}]^{\ell-1}$ only for $\ell=1$. It means that in the
notations of Section 2.4 we should consider the partitions
$(k_1,\dots,k_{\ell})$, $k=k_1+\dots+k_\ell$ only with $\ell=1$.
Thus, we have proved
\subsubsection{}
\begin{lemma}
Let $\F_k$, $k\ge 1$ be the Taylor components of the $A_\infty$-morphism
$\F\colon\Dif^\mb(\E)\to
\End(K_0^\mb)$. Then the supertrace $\Upsilon=\F_{\Hoch}^*(\str_{can})$ can be
described as follows: its value $\Upsilon(a)$ for
$a\in[\Dif^\mb(\E)^{\otimes k}]^{k-1}$ is
\begin{equation}\label{eqsem33}
\Upsilon(a)=\sum_{s=0}^{k-1}\str_{can}(\F_k(C^s(a)))
\end{equation}
where $C$ is the cyclic shift defined on $\cA^{\mb\otimes k}$ as
\begin{equation}\label{eqsem34}
C(a_1\otimes\dots\otimes a_k)=(-1)^{(\deg a_k+1)(\deg a_1+\dots+\deg
a_{k-1}+k-1)}a_k\otimes a_1\otimes\dots\otimes a_{k-1}
\end{equation}
\begin{proof}
It follows from Lemma 2.4.
\end{proof}
\end{lemma}

\subsubsection{Topological Quantum Mechanics and Hodge theory}
Here we specify the formula (\ref{eqsem33}) for the trace even more, taking
in the account a Kahler metric on the manifold $X$ and an Hermitian metric in
the bundle $E$, and the corresponding Hodge theory. For simplicity we
consider here only the case $E=\mathbb{C}$.

Now let $\cA^\mb=\Dif^\mb(\cO)$ and $M^\mb=\Gamma_X(\Dolb^\mb(X,\cO))$, $Q=\bar{\partial}$.
Suppose a Kahler metric on $X$ is chosen. Denote by $\Delta$ the
$\bar{\partial}$-Laplacian. As before, we suppose that $X$ is compact.
Set
\begin{equation}\label{eqsem35}
\begin{aligned}
\ &M_0^\mb=\K \Delta={harmonic\ \  forms}\\
&M^\mb_1=\I (\bar{\partial})\oplus\I (\bar{\partial}^*)
\end{aligned}
\end{equation}
According to the classical Hodge theory, $M^\mb=M^\mb_0\oplus M^\mb_1$,
the differential on $M_0^\mb$ is 0, and $M_1^\mb$ is acyclic.
Choose the homotopy $\kappa\colon M^\mb\to M^\mb_1[-1]$ as
$\bar{\partial}^*$. Then $\{Q,\kappa\}=\Delta$, the
$\bar{\partial}$-Laplacian. So we replace $\Pi_{M_1}$ by the Laplacian, that
is, we are in the situation of the Remark 2.3.1.
Only what we need to check is the equation (iii) in this Remark, namely,
\begin{equation}\label{eqsem36}
\exp(-t\Delta)(m_1)\rightarrow 0 \ \  {when} \ \ t\rightarrow +\infty\ \
{for}\ \
\forall m_1\in M^\mb_1
\end{equation}
To prove (\ref{eqsem35}) notice that the Laplacian $\Delta$ on any compact
Kahler manifold has a {\it descrete} non-negative spectrum. The space $M_0^\mb$ is
exactly the 0 eigenspace. Let $\lambda_1$ be the smallest positive
eigenvalue of $\Delta$. Then
\begin{equation}\label{eqsem37}
\|\exp(-t\Delta)(m_1)\|\le \exp(-t\lambda_1)\cdot\|m_1\|
\end{equation}
It is clear that for a fixed $m_1$, the r.h.s. of (\ref{eqsem37}) tends to 0
when $t\rightarrow +\infty$.

In the next Sections we extend the super-trace $\Upsilon$ to some
completion of the Hochschild chain complex $\Hoch_\mb(\Dif^\mb(\E))$. The
difficulty here is that for this completed complex we have {\it not} any map
to the Hochschild complex of the algebra
$\End K_0^\mb$. Therefore, we can not anymore construct our supertrace
$\widehat{\str}$ by the pull-back of some trace on the algebra of endomorphisms. The
only what remains is to look at the formula (\ref{eqsem33}) and to try to
extend {\it it} to the completed complex. Then this extension will be a
trace (a linear functional on 0-th cohomology of the completed Hochschild
complex) automatically.
\bigskip
\subsection{The completed Hochschild complex
$\widehat{\Hoch_\mb(\Dif^\mb(\E))}$}
We are going to describe here a compactification of the Hochschild complex
$\Hoch_\mb(\Dif^\mb(\E))$ of the algebra $\Dif^\mb(\E)$ which gives the
"right" homology. We start with a simpler question: Let $M$ be a
$C^\infty$-manifold (of dimension $n$, or a dg manifold), what is the right
definition of the Hochschild homology of the algebra $A=C^\infty(M)$? By
the right definition we mean such a definition that the
Hochschild-Kostant-Rosenberg theorem holds, namely, the homology $\HH_i(A)$
are isomorphic to $i-th$ smooth differential forms on $M$, $\Omega^i(M)$.
This is well known that there are several such definitions of
$\Hoch_\mb(A)$. Recall here all of them:
\begin{itemize}
\item[1.] $A^{\otimes k}=C^\infty(M^k)$
\item[2.] $A^{\otimes k}={germs}_{\Delta}C^\infty(M^k)$
\item[3.] $A^{\otimes k}={jets}_{\Delta}C^\infty(M^k)$
\end{itemize}
(here $\Delta\colon M\hookrightarrow M^k$ is the diagonal).
In particular, the naive algebraic definition $A^{\otimes k}$ as the usual
algebraic tensor power of $A$ does not work in the $C^\infty$-case. It works
when $M=V$ is a vector space, and $A$ is the algebra of {\it polynomial}
functions on $V$.

We have
the following
\begin{lemma1}
In all three definitions above the Hochschild homology
$\HH_i(A)=\Omega^i(M)$,
the $C^\infty$ differential forms.
\qed
\end{lemma1}

In the present paper we use the first definition from the list above
basically because we want to define a completion. Notice that from the point
of view of proofs the second and the third definitions are simpler because
we can use the sheaf theory in the cohomology computations.

Now we pass to differential operators in the $C^\infty$ sense (in trivial
bundle).
\begin{lemma2}
Let $M$ be a $C^\infty$-manifold, $\dim_{\mathbb{R}}M=n$, and let $A=\Dif(M)$ be the algebra of
$C^\infty$ differential operators on $M$. Define $A^{\otimes k}=\Dif(M^k)$.
Then the Hochschild homology $HH_{-i}(\Dif(M))=\CH^{2n-i}(M)$.
\begin{proof}
One way to prove this Lemma is to localize Lemma 1 in the sense of the
Tsygan formality (see [Sh]). Another way is to prove somehow (maybe again
using the Tsygan formality) the local statement for $M=V$ and then to use
the simplicial methods. Here the main difficulty is that we should extend
the Eilenberg-Zilber theorem for the completed tensor products. Anyway, if
we fix the isomorphism locally, we get a canonical global isomorphism.
\end{proof}
\end{lemma2}
\begin{corollary}
Let $M$ be a $C^\infty$-manifold, $A=\gl_N(\Dif(M))$, and $A^{\otimes k}=
\gl_N^{\otimes k}\otimes \Dif(M^k)$. Then again $\HH_{-i}(A)=\CH^{2n-i}(M)$.
\begin{proof}
Use the Kunneth formula for the Hochschild homology.
\end{proof}
\end{corollary}
Our main result here is
\begin{theorem*}
Let $X$ be a complex manifold, $\dim_{\mathbb{C}}M=n$, $E$ be a holomorphic vector bundle over it,
$\E$ be the corresponding sheaf of holomorphic sections. Define now the $k$-th
tensor power $\Dif^\mb(\E)^{\otimes k}$ as the $\Dif(\E^{\boxtimes
k})\otimes_{C^\infty(X)^{\otimes k}}C^\infty(T^{0,1}[1]X)^{\otimes k}$ where $\Dif$ stands for $C^\infty$-differential operators.
Then the corresponding chain Hochschild complex has the
homology $\widehat{\HH_{-i}}(\Dif^\mb(\E))=\CH^{2n-i}(X)$. There exists a
canonical (depending on $E$) isomorphism
$\vartheta_E\colon\widehat{\HH_{-i}}(\Dif^\mb(\E))\to\CH^{2n-i}(X)$.
\begin{proof}
We use simplicial methods and the Eilenberg-Zilber theorem (actually, a
"completed" version of it). For the local computation we use the Corollary
above. It gives us automatically the map $\vartheta_E$.
\end{proof}
\end{theorem*}
\begin{remark}
We used here the first definition (in the listing at the beginning of
Section 3.2) of the completed Hochschild homology basically because of
aesthetic reasons. We mean that this definition more than others two is
associated with a {\it completion}. If we use the third definition we could
escape the simplicial methods and a completion of the Eilenberg-Zilber
theorem. Instead of it, we can use elementary sheaf methods. We mean that in
this (the third) definition we can consider $A^{\otimes k}$ as a sheaf on
$X$ (not on $X\times\dots\times X$). Then we can simply globalize the local
computation, for which we need to have the Hochschild-Kostant-Rosenberg
theorem.
\end{remark}
\subsection{The Integral conjecture}
\subsubsection{The trace $\widehat{\str}$}
First of all, we prove the following
\begin{theorem*}
The formula (\ref{eqsem33}) (in the explicit homotopy of Section 3.1.2) for
the super-trace $\Upsilon=\F^*_{\Hoch}(\str_{can})$ on the algebra
$\Dif^\mb(\E)$ (which is nonzero only on $[\Dif^\mb(\E)^{\otimes k}]^{k-1}$,
$k\ge 1$) can be continued to the completed tensor power $\widehat{[\Dif^\mb(\E)^{\otimes
k}]}^{k-1}$ defined in Section 3.2, and then it defines a trace $\widehat{\str}$ on the completed
Hochschild chain complex of the algebra $\Dif^\mb(\E)$.
\begin{proof}
Let $D=D_k\otimes\dots\otimes D_1\in [\Dif^\mb(\E)^{\otimes k}]^{k-1}$ be a
typical indecomposable chain before the completion. We can rewrite the
formula (\ref{eqsem33}) for $\Upsilon(D)$ as follows. First, we can
associate with $D$ canonically an object on $X^{\times k}$, namely, an
element of $\Dif(\E^{\boxtimes
k})\otimes_{C^\infty(X)^{\otimes k}}C^\infty(T^{0,1}[1]X)^{\otimes k}$.
Denote this element by $\cD$. We want to rewrite our formula for
$\Upsilon(\cD)$ to make it sense for an arbitrary ("completed") $\cD\in
\Dif(\E^{\boxtimes
k})\otimes_{C^\infty(X)^{\otimes k}}C^\infty(T^{0,1}[1]X)^{\otimes k}$, not
only for $\cD$ equal to the tensor product of $k$ factors. It is sufficiently
to do it for $\F_k(\cD)$ because of the formula (\ref{eqsem33}). As in
Section 3.1.2, we introduce a Kahler metric on $X$ and consider the
$\bar{\partial}$-Laplacian $\Delta$, and the corresponding Hodge theory.
Then we can rewrite formula (\ref{eqr1.12}) as
\begin{multline}\label{eqr1.55}
\Omega_{\cD}=\Pi_{K^\mb_0}\circ D_k\circ \exp[-d(\tau_{k}-\tau_{k-1})\bar{\partial}^*-
(\tau_{k}-\tau_{k-1})\Delta]\circ D_{k-1}\circ \\
\exp[-d(\tau_{k-1}-\tau_{k-2})\bar{\partial}^*-
(\tau_{k-1}-\tau_{k-2})\Delta]\circ D_{k-2}\circ \\
\dots\circ\exp[-d(\tau_2-\tau_1)\bar{\partial}^*-
(\tau_2-\tau_1)\Delta]\circ D_1\circ \Pi_{K^\mb_0}\in \End^\mb(K^\mb_0)\otimes
\Omega^\mb(C_{k+1})
\end{multline}
Now introduce
\begin{multline}\label{eqsem1.57}
\Phi=\exp[-d(\tau_{k}-\tau_{k-1})\bar{\partial}^*-
(\tau_{k}-\tau_{k-1})\Delta]_{k-1}\circ\dots\\
\circ\exp[-d(\tau_2-\tau_1)\bar{\partial}^*-
(\tau_2-\tau_1)\Delta]_1
\end{multline}
(we consider $\Phi$ as an object on $X\times X\times\dots\times X$, the lower index denotes the variable at
which the operator acts, and $\circ$ denotes the composition of operators).
It is clear that $\Phi$ is {\it not} a differential operator. Nevertheless,
the action of it on the $k$-th tensor power of the Dolbeault complex $K^{\mb\otimes k-1}=\Gamma_X(\Dolb^\mb(X,\E))^{\otimes k-1}$
is well-defined. Moreover, as operators acting on different variables
commute, we can rewrite (\ref{eqr1.55}) as
\begin{equation}\label{eq_super}
\Omega_{\cD}=\Pi_{K^\mb_0}\circ m(\Phi\circ (D_k\otimes\dots\otimes D_1))\circ \Pi_{K^\mb_0}
\end{equation}
where $m$ is the restriction to the diagonal $X\hookrightarrow
X\times\dots\times X$, or, in other words, the product.
The last formula clearly makes sense for an arbitrary $\cD\in\Dif(\E^{\boxtimes
k})\otimes_{C^\infty(X)^{\otimes k}}C^\infty(T^{0,1}[1]X)^{\otimes k}$, and
defines a trace on the completed Hochschild complex.
\end{proof}
\end{theorem*}
\begin{remark}
We used in this proof that $\Delta$ is a positive operator with a descrete
spectrum.
\end{remark}
\subsubsection{The integral conjecture and the Riemann-Roch-Hirzebruch
theorem}
In Theorem 3.3.1 we constructed a trace $\Upsilon\colon
\mathbb{H}^0(\widehat{\Hoch_\mb}(\Dif^\mb(\E))\to\mathbb{C}$, and in Theorem
3.2 we proved that
$\mathbb{H}^0(\widehat{\Hoch_\mb}(\Dif^\mb(\E))\simeq\CH^{2n}(X)$ (the last
isomorphism is $\vartheta_E$). These two results give us a map
\begin{equation}\label{eqsem58}
\Im\colon\CH^{2n}(X)\to\mathbb{C}
\end{equation}
\begin{intconj}
The map $\Im$ is the integral over the fundamental cycle of the (compact)
manifold $X$.
\end{intconj}
This statement looks quite mysterious: we constructed the integral starting
from absolutely different things. We can prove this Conjecture only in the
case when $\chi(\E)\ne 0$, and using the RRH theorem. Moreover, we need the
following
\begin{lemma}
Consider the map
$\vartheta_E\colon\widehat{\HH_{0}}(\Dif^\mb(\E))\to\CH^{2n}(X)$ from Theorem
3.2. Let $\Id$ be the identity differential operator. Then $\vartheta(\Id)=
[Td(T_X)\cdot ch(E)]_{2n}$.
\begin{proof}
This statement is proven, in a form, in many places. See in particular [FT],
[NT1], [FFSh].
\end{proof}
\end{lemma}
Now the following statement is almost a tautology:
\begin{theorem*}
The Integral Conjecture is true when $\chi(\E)\ne 0$.
\begin{proof}
It follows from Lemma above and from the RRH theorem that if $\chi(\E)\ne 0$
the identity differential operator $\Id$ does {\it not} belong to the kernel
of the map $\vartheta_E$. Then the image $\vartheta_E(\Id)$ is a non-zero
element in the 1-dimensional space $\CH^{2n}(X)$. By the construction,
\begin{equation}\label{eqsem59}
\str(\cD)=\Im (\vartheta_E(\cD))
\end{equation}
For any holomorphic differential operator $\cD$ in $\E$.
Applying this to $\cD=\Id$ and using the RRH theorem obtain
$\Im (\vartheta_E(\Id))=\int_X(\vartheta_E(\Id))$. Then we conclude that
$\Im=\int_X$ because $\vartheta_E(\Id)\ne 0$ and because the space
$\CH^{2n}(X)$ is 1-dimensional.
\end{proof}
\end{theorem*}
We have
\begin{corollary}
Suppose that $\chi(\E)\ne 0$. Then for any holomorphic differential operator $\cD\in\Dif(\E)$ one has
\begin{equation}\label{eqsem60}
\str(\cD)=\int_X(\vartheta_E(\cD))
\end{equation}
\qed
\end{corollary}
The last formula can be rewritten explicitely using the Hochschild cocycle
from [FFSh] and formal geometry.

The reader could notice that here the Integral Conjecture has the primary
interest for us, while the RRH theorem and formula (\ref{eqsem60}) play an
auxiliary role. It would be very interesting to find a direct approach to
the Integral Conjecture, and, in particular, to prove it without the
assumption $\chi(\E)\ne 0$.

In the next Section we consider the analog of the Integral Conjecture in the
case of cyclic homology. Here we have not any RRH theorem, and even the
complete formulation of the conjecture still remains to be open.
\section{The case of cyclic homology}
Consider the cyclic homology instead of the Hochschild homology in all
constructions above.

The $A_\infty$-morphism $\F\colon\Dif^\mb(\E)\to K^\mb_0$ induces the map
$\F_{\Cycl}\colon\Cycl_\mb(\Dif^\mb(\E))\to\Cycl_\mb(\End K_0)$ (see Section
2.4).

We have the following
\begin{lemma1}
Let $\mathbf{A}^\mb$ be the algebra of endomorphisms of a finite-dimensional
graded (with zero differential) vector space $V^\mb$. Then the cyclic
homology $\mathbb{H}^i(\Cycl_\mb(\mathbf{A}^\mb))$ is equal to $\mathbb{C}$
for an even not-positive $i$ and is equal to 0 otherwise. One has a
canonical functional
$\Tr_{2i+1}\colon\mathbb{H}^{-2i}(\Cycl_\mb(\mathbf{A}^\mb))\to\mathbb{C}$
given by the formula
\begin{equation}\label{eqsem70}
\Tr_{2i+1}(A_1\otimes A_2\otimes\dots\otimes
A_{2i+1})=\Tr\left(\sum_{\sigma\in\Sigma_{2i+1}}(-1)^{\natural(\sigma)}A_{\sigma(1)}\cdot\dots\cdot
A_{\sigma(2i+1)}\right)
\end{equation}
\qed
\end{lemma1}
The advantage comparably with the case of the Hochschild homology is that
here we have "many" traces, namely, a trace in each not-positive even degree,
for the algebra $\mathbf{A}^\mb$ of endomorphisms, while in the Hochschild
case we have such a trace only in degree 0.

We can take the pull-back $\Upsilon_{2i}=\F_{\Cycl}^*(\Tr_{2i+1})$. This is
a "higher trace", that is, a linear functional
$$
\Upsilon_{2i}\colon\mathbb{H}^{-2i}(\Cycl_\mb(\Dif^\mb(\E)))\to\mathbb{C}
$$
Now it is interesting to compute the cyclic homology
$\mathbb{H}^{\mb}(\Cycl_\mb(\Dif^\mb(\E)))$. As in the Hochschild case, the
answer is not interesting before the completion.
There exists a completion of the cyclic complex $\Cycl_\mb(\Dif^\mb(\E))$
analogous to the completion of the Hochschild complex described in Section
3.2. Denote the completed cyclic complex by
$\widehat{\Cycl_\mb}(\Dif^\mb(\E))$.
We have
\begin{lemma2}
The cohomology $\mathbb{H}^{-i}(\widehat{\Cycl_\mb}(\Dif^\mb(\E)))$ is equal
$\CH^{2n-i}(X)\oplus\CH^{2n-i+2}(X)\oplus\CH^{2n-i+4}(X)\oplus\dots$.
\qed
\end{lemma2}
In particular,
$\mathbb{H}^{0}(\widehat{\Cycl_\mb}(\Dif^\mb(\E)))\simeq\CH^{2n}(X)$ and
$\mathbb{H}^{-2}(\widehat{\Cycl_\mb}(\Dif^\mb(\E)))\simeq\CH^{2n}(X)\oplus\CH^{2n-2}(X)$.

We get a map
\begin{equation}\label{eqsem70}
\Im_{-2}\colon\CH^{2n}(X)\oplus\CH^{2n-2}(X)\to\mathbb{C}
\end{equation}

What is this map? The construction depends on the bundle $E$ and therefore
the map $\Im_{-2}$ also could depend. Conjecturally, its component
$\Im_{-2}^{2n}\colon\CH^{2n}(X)\to\mathbb{C}$ is the integral (at least, up
to a constant depending only on the dimensions). But we have no tools in the
moment to describe the map
$\Im_{-2}^{2n-2}\colon\CH^{2n-2}(X)\to\mathbb{C}$. It could be 0 but could be not. The only what we can say
is the following:

Let $\cD_1\otimes\cD_2\otimes\cD_3$ be a cyclic cycle, $\cD_i$ be
holomorphic differential operators in $\E$.
Then
\begin{equation}\label{eqsem71}
\str\left(\sum_{\sigma\in\Sigma_3}(-1)^{\natural
\sigma}\cD_{\sigma(1)}\circ\cD_{\sigma(2)}\cD_{\sigma(3)}\right)=\Im_{-2}^{2n-2}([\cD_1\otimes\cD_2\otimes\cD_3])
\end{equation}
where $[\cD_1\otimes\cD_2\otimes\cD_3]\in\CH^{2n-2}(X)$ be the corresponding
class. This formula follows directly from the definitions.

We can compute $[\cD_1\otimes\cD_2\otimes\cD_3]$ in one particular case,
namely, we can compute $[\Id\otimes \Id\otimes \Id]$. It was proven in Nest-Tsygan
papers that $[\Id\otimes\Id\otimes\Id]=[Td(T_{X})\cdot ch(E)]_{2n-2}$ and
$[\Id]=[Td(T_{X})\cdot ch(E)]_{2n}$.
Therefore, if one belives that $\Im_{-2}^{2n}$ is the integral,
$\Im_{-2}^{2n-2}([Td(T_X)\cdot ch(E)]_{2n-2})$ should be 0 by the RRH theorem. In general, we
expect that $\Im_{-2}^{2n-2}$ is the integral over some combination in
codimension 2 of the Poincare dual to the Chern classes of the bundle $E$,
and of the natural bundles associated with $X$.

It would be very interesting to understand better this topic.
\subsection*{Acknowledgements}
The authors are grateful to Giovanni Felder for many discussions,
suggestions and improvements. The work of A.L. is partially supported by
Russian Federal Program 40.052.1.1.1112 and by the grants: Volkswagen
Stiftung, RFFI 01-01-00548 and NSh-1999/2003.2
\newpage

\bigskip
Independent University of Moscow,\\
11 Bolshoj Vlas'evskij pereulok, 121002 Moscow, RUSSIA\\
{\it e-mail}: {\tt feigin@mccme.ru}\\
\\
Moscow Institute for Theoretical and Experimental Physics (ITEP),\\
25 Bolshaja Cheremushkinskaja ulica, Moscow RUSSIA\\
{\it e-mail}: {\tt losev@mail.itep.ru}\\
\\
Department of Mathematics, ETH Zurich,\\
CH-8092 Zurich SWITZERLAND\\
{\it e-mail}: {\tt borya@math.ethz.ch}\\
\end{document}